\documentclass{amsart}
\usepackage{amsmath}%
\usepackage{amsfonts}%
\usepackage{amssymb}%
\usepackage{graphicx}
\newtheorem{theorem}{Theorem}
\theoremstyle{plain}

\newtheorem{corollary}{Corollary}

\newtheorem{lemma}{Lemma}

\newtheorem{proposition}{Proposition}

\numberwithin{equation}{section}

\def\cH{\mathcal{H}}

\def\IC{\mathbb{C}}

\def\IP{\mathbb{P}}

\begin{document}
\title[Complexity of Bezout's Theorem VI]{Complexity of Bezout's Theorem VI: Geodesics in the Condition (number) metric}
\author{Michael Shub}
\address {Department of Mathematics, University of Toronto, 40 St. George Street, Toronto, Ontario, M5S 2E4,
Canada}
\email{michael.shub@utoronto.ca}
\urladdr{http://www.math.toronto.edu/shub}
\thanks{ This work was partly supported by an NSERC Discovery Grant.}

\date{February 1,2007}

\begin{abstract}
We introduce a new complexity measure of a path of (problems,
solutions) pairs in terms of the length of the path in the condition
metric which we define in the article. The measure gives an upper
bound for the number of Newton steps sufficient to approximate the
path discretely starting from one end and thus produce an
approximate zero for the endpoint. This motivates the study of short
paths or geodesics in the condition metric.
\end{abstract}
\maketitle

\section{Introduction}

\noindent In a series of papers  we have studied the complexity of
solving systems of homogeneous polynomial equations by applying
Newton's method to a homotopy of (system,solution) pairs
\cite{BezI}, \cite{BezII},\cite{BezIII}, \cite{BezIV}, \cite{BezV}.
The latest word in this direction is \cite{Beltran-Pardo I},\cite{
Beltran-Pardo II}. A key ingredient is an estimate of the number of
Newton steps by the maximum condition number  along the path
multiplied by the length of the path of solutions. The main result
of this paper is to show that the maximum condition number times the
length may be replaced by the integral of the condition number times
the length of the tangent vector to the path, Theorem \ref
{MainTheorem}. The result suggests using the condition number to
define a Riemannian metric on the solution variety and to study the
geodesics of this metric. Finding a geodesic is in itself not easy.
So we do not immediately obtain a practical algorithm. Rather the
study may help to understand in some systematic fashion the geometry
of homotopy algorithms, especially those that attempt to avoid
ill-conditioned problems. We note that recentering algorithms in
linear programming theory  may be seen as adaptively avoiding
ill-conditioning. \cite{Nesterov-Todd} compares the central path of
linear programming theory to geodesics in an appropriate metric. In
\cite{Beltran-Shub} we begin studying distances in the condition
metric.

\section{Definitions and Theorems}

\noindent We begin by recalling the context. For every positive
integer $l\in{\mathbb N}$, let $H_l\subseteq {\mathbb
C}[X_0,\ldots,X_n]$ be the vector space of all homogeneous
polynomials of degree $l$.  For $(d):=(d_1,\ldots,d_n)\in{\mathbb
N}^n$, let ${\mathcal H}_{(d)}:=\prod_{i=1}^n H_{d_i}$ be the set of
all systems $f:=(f_1,\ldots,f_n)$ of homogeneous polynomials of
respective degrees $\deg(f_i)=d_i,~1\leq i\leq n$. So $f:
\mathbb{C}^{n+1}\rightarrow \mathbb{C}^n.$ We denote by
$D:=\max\{d_i:1\leq i\leq n\}$ the maximum of the degrees.

The solution variety $\hat{V} \subset \mathcal {H}_{(d)}\times
(\mathbb{C}^{n+1}-\{0\}) $ is the set of points $\{(f,x)|f(x)=0\}.$
Since the equations are homogeneous, for all $\lambda_1,\lambda_2
\in \mathbb{C}-\{0\}$, $\lambda_1f(\lambda_2 x) =0$ if and only if
$f(x)=0.$ So $\hat{V}$ defines a variety $V\subset
\mathbb{P}({\mathcal H}_{(d)})\times \mathbb{P}(\mathbb{C}^{n+1})$
 where $\mathbb{P}({\mathcal H}_{(d)})$ and
 $\mathbb{P}(\mathbb{C}^{n+1})$ are the projective spaces
 corresponding to ${\mathcal H}_{(d)}$ and $\mathbb{C}^{n+1}$
 respectively. $\hat{V}$ and $V$ are smooth. We speak
 interchangeably of a path $(f_t, \zeta_t)$ in $\hat{V}$ and its
 projection $(f_t, \zeta_t)$ in $V.$ Most quantities we define are
 defined on $\hat{V}$ but are constant on equivalence classes so are
 defined on $V.$

 For $(g,x)\in \mathbb{P}({\mathcal H}_{(d)})\times
 \mathbb{P}(\mathbb{C}^{n+1})$ let
 $$ \mu_{norm}(g,x) = ||g||
 ||Dg(x)|_{N_x}^{-1}\Delta(d^{\frac{1}{2}}_i||x||^{d_i-1})||$$

 where $||g||$ is the unitarily invariant norm defined by the
 unitarily invariant Hermitian structure on ${\mathcal H}_{(d)}$ considered in
 \cite{BezI} and sometimes called the Bombieri-Weyl or Kostlan
 Hermitian structure, $||x||$ is the standard norm in
 $\mathbb{C}^{n+1}$, $N_x$ is the Hermitian complement of $x$, and
 $\Delta(a_i)$ for $a_i \in \mathbb{C}, i=1 \dots n$ is the $n
 \times  n$ diagonal matrix with $i th$ diagonal entry $a_i.$ If
 $Dg(x)|_{N_x}^{-1}$ does not exist we take $\mu_{norm} = \infty.$
 $\mu_{norm}$, also called $\mu_{proj}$ in some of our papers, is a
 normalized version of the condition number which we usually denote
 by $\mu.$ We have also used various notions of distance in
 projective space. In this paper we use only the Riemannian distance inherited from the
 Hermitian structure on the vector space, i.e. the angle. In
 previous papers we have paid careful attention to the constants.
 In this paper we are more cavalier.
  We begin with an analysis of how the normalized condition
 number varies.

\begin{theorem}\label{Theorem 1}
Given $\epsilon > 0$ there is a constant $C > 0$ such that if

$$g \epsilon \IP (\cH_{(d)})  \ \rm and \ \zeta, \eta \epsilon \IP(\IC^{n+1})$$

with

$$d (\zeta, \eta) < \frac{C}{D^{3/2} \mu_{norm} (g, \zeta)}$$

Then $\mu_{norm} (g ,\eta) \le (1 + \epsilon ) \mu_{norm} (g,
\zeta).$
\end{theorem}

\medskip

\begin{proof}
Use proposition 2.3 of \cite{BezIV}.

Note that $u \le D^{3/2} \mu_{norm} (g, \eta) \cdot d (\eta, \zeta)
\le C$.  $r_0 \sim d (\eta, \zeta)$, the $\eta(f,x)$ in the
definitions of $K$ is $\le 1$ and
$\bigl(\frac{\|\eta\|}{\|\zeta\|}\bigr)^{D-1} \sim (1 + d (\zeta,
\eta))^{D-1} < e^c$.
\end{proof}

\bigskip

\begin{theorem}\label{Theorem 2}
Given $\epsilon > 0$ there is a constant $C > 0$ such that if

$$f, g \epsilon \IP (\cH_{(d)})\  \rm and \
\zeta \epsilon \IP (\IC^{n+1}) \ \rm $$

and

$$d (f, g) < \frac{C}{D^{1/2} \mu_{norm} (f, \zeta)},$$

then

$$\mu_{norm} (g, \zeta) \le (1 + \epsilon) \mu_{norm} (f, \zeta)$$
\end{theorem}

\begin{proof}
By Proposition 5b) of \cite{BezI} SectionI-3.

\begin{align*}
\mu_{norm} (g, \zeta) &\le \frac{\mu_{norm} (f, \zeta) (1 + d (f, g))}{1-D^{1/2} d(f, g) \mu_{norm}
(f, \zeta)} \\
&\le \frac{\mu_{norm} (f, \zeta) \bigl(1 + \frac{C}{D^{1/2} \mu_{norm} (f, \zeta)} \bigr)}{1-C}
\end{align*}

Recall that $\mu_{norm} (f, \zeta) \ge 1$.
\end{proof}

\bigskip

\begin{theorem}\label{Theorem 3}
Given $\epsilon > 0$ there is a constant $C > 0$ such that if $$f, g
\epsilon \IP (\cH_{(d}))\ \rm and \  \zeta, \eta \epsilon \IP
(\IC^{n+1})\ \rm $$ and

\begin{align*}
d (f, g) &< \frac{C}{D^{1/2} \mu_{norm} (f, \zeta)}\\
d (\zeta, \eta) &< \frac{C}{D^{3/2} \mu_{norm} (f, \zeta)},
\end{align*}

then

$$\frac{1}{1+\epsilon} \mu_{norm} (g, \eta) \le \mu_{norm} (f, \zeta) \le
(1+\epsilon) \mu_{norm} (g, \eta).$$
\end{theorem}

\medskip
\begin{proof}
Apply Theorems 1 and 2 to prove the left hand inequality.  Then
given the left hand inequality apply Theorem 1 and 2 again to prove
the right hand inequality applying the theorems with $g, \eta$ in
place of $f, \zeta$.  Adjust $C$ and $\epsilon$ as necessary.
\end{proof}

\bigskip

The next proposition is useful for our Main Theorem \ref
{MainTheorem}.

\begin{proposition}
Given $\epsilon > 0$ there is a $C > 0$ with the following property:

\smallskip
Let $(f_t, \zeta_t)$ be a $C^1$ path in $V$ for $t_0 \leq t \leq
t_1$. Define $S_0 = t_0$ and $S_i$ to be the first value of $t \le
t_1$ such that

$$\int\limits_{S_{i-1}}^{S_i}
(\| \dot f_t \| + \| \dot \zeta_{t} \|) d_t = \frac{C}{D^{3/2}
\mu_{norm} (f_{S_{i-1}}, \zeta_{S_{i-1}})} \ {\rm or} \ t_1.$$

Then $S_k = t_1$ for
$$k \le max(1, \frac{(1+\epsilon)}{C} D^{3/2} \int\limits_{t_0}^{t_1} \mu_{norm} (f_t,
\zeta_t) (\| \dot f_{t} \| + \| \dot \zeta_{t} \|) dt)$$

and $\mu_{norm} (f_{t_1}, \zeta_{t_1}) \le (1 + \epsilon)^k
\mu_{norm} (f_{t_0}, \zeta_{t_0})$
\end{proposition}

\medskip

\begin{proof}
\begin{align*}
\int\limits_{S_{i-1}}^{S_i} \mu_{norm} (f, \zeta_t) (\| \dot f_{t}
\| + \| \dot \zeta_{t} \|) dt &\ge \frac{1}{(1+\epsilon)}
\int\limits_{S_{i-1}}^{S_i} \mu_{norm} (f_{S_{i-1}} ,
\zeta_{S_{i-1}})
\ (\| \dot f_{t} \| + \| \dot \zeta_{t} \| ) dt \\
&\ge \frac{{\frac{1}{1+\epsilon}} C}{D^{3/2}} \ {\rm if} \ S_i < t_1.
\end{align*}

Consequently

$$\int\limits_{S_0}^{S_k} \mu_{norm} (f_t, \zeta_t) \ (\| \dot f_{t} \| + \| \dot \zeta_{t} \|) dt$$

$$\ge \frac{(k-1) C}{(1 + \epsilon) D^{3/2}} \ {\rm and}$$

$$(k-1) \frac{ C}{(1 + \epsilon)D^{3/2}} \le \quad \int\limits_{t_0}^{t_1} \mu_{norm} (f_{t_1} \zeta_{t})
\ (\| \dot f_{t} \| + \| \dot \zeta_{t} \|) dt.$$

$${\rm so} \ k-1 \le \frac{(1+\epsilon) D^{3/2}}{C} \int\limits_{t_0}^{t_1} \mu_{norm}
(f_{t_1} \zeta_{t}) \ (\| \dot f_{t} \| + \| \dot \zeta_{t} \|) dt.$$
\end{proof}

Since we are working in the metric $d$ we require an approximate
zero theorem in this metric. First we prove a lemma. Recall the
following quadratic polynomial $$\psi(u)=1-4u+2u^2$$ and the
definition of the projective Newton iteration $$N_f(x)=
x-(Df(x)|_{x^\perp})^{-1}f(x).$$ Here $x^\perp$ is the Hermitian
complement to $x.$ $N_f:\IP(\IC^{n+1})\rightarrow \IP(\IC^{n+1})$
except that it fails to be defined where $(Df(x)|_{x^\perp})^{-1}$
does not exist. If $f(\zeta)=0$ and $d (N^k_f (x), \zeta) \le
\frac{1}{2^{2^k-1}} d (x, \zeta)$ for all positive integers $k$,
then $x$ is called an $\emph{approximate zero}$ of $f$ with
associated zero $\zeta.$

\medskip
\begin{lemma}
Let $u < \frac{3 - \sqrt{7}}{4}$.  Let $f \epsilon P (\cH_{(d)})$
and $\zeta \epsilon P  (\IC^{n+1})$ with $f (\zeta) = 0$. If $d (x,
\zeta) \le \frac{u}{D^{3/2} \mu_{norm} (f, \zeta)},$ then

$$d(N_f (x), \zeta) \le \frac{4u}{\psi (2 u)} d (x, \zeta).$$
\end{lemma}

\medskip

\begin{proof}
In the range of angles under consideration $\rm Tan \ d(x, \zeta) \le 2 d (x, \zeta)$.
So apply Lemma 1 of P263 of \cite{BCSS} to conclude that

\begin{align*}
d (N_f (x), \zeta) &\le {\rm Tan} \ d (N_f (x), \zeta) \le \frac{2
u}{\psi (2 u)} {\rm Tan} \
d(x, \zeta)\\
&\le \frac{4u}{\psi (2u)} d(x, \zeta).
\end{align*}
\end{proof}

Now let $u_0$ solve the equation

$$\frac{4 u_0}{\psi (2 u_0)} = \frac{1}{2} \quad \rm or \ u_0 =
\frac{16 - \sqrt{232}}{16} \sim 0.048$$

\medskip

\begin{theorem}\label{Theorem 5} (Approximate zero Theorem)

\medskip
Let $f \epsilon P (\cH_{(d)}), \ \zeta \epsilon P \ (\IC^{n+1})$
with $f(\zeta) = 0$.  If $d (x, \zeta) < \frac{u_0}{D^{3/2}
\mu_{norm} (f, \zeta)}$ Then $d (N^k_f (x), \zeta) \le
\frac{1}{2^{2^k-1}} d (x, \zeta).$
\end{theorem}

\medskip

\begin{proof} By induction

\medskip
Suppose

\begin{align*}
d(N^k_f (x), \zeta) &\le \bigl(\frac{4 u_0}{\psi (2 u_0)}\bigr)^{2^k-1} d(x, \zeta)\\
&\le \frac{\bigl(\frac{4 u_0}{\psi (2 u_0)}\bigr)^{2^k-1}
u_0}{D^{3/2} \mu_{norm} (f, \zeta).}
\end{align*}

Let $u_1 = \bigl(\frac{4 u_0}{\psi (2 u_0)}\bigr)^{2^k-1} u_0$. Note
$u_1 \le u_0$ and so $\psi (2 u_1) \ge \psi (2 u_0)$.  By the lemma

\begin{align*}
d (N^{k+1}_f(x), \zeta) &\le
\frac{4 u_1}{\psi (2 u_1)} d (N^k_f (x), \zeta)\\
&\le \frac{4 u_1}{\psi (2 u_0)} \cdot d (N^k_f (x), \zeta) \\
&\le 4 \frac{\bigl(\frac{4 u_0}{\psi (2 u_0)}\bigr)^{2^{k}-1}
u_0}{\psi (2 u_0)}
\cdot \bigl(\frac{4 u_0}{\psi (2 u_0)}\bigr)^{2^{k}-1} d (x, \zeta) \\
&= \bigl(\frac{4 u_0}{\psi (2 u_0)}\bigr)^{2^{k+1}-1} d (x, \zeta)\\
\end{align*}
\end{proof}

\bigskip
Let $(f_t, \zeta_t)$ be a (piecewise) $C^1$ path in $V, \
\frac{d}{dt} (f_t, \zeta_t)$ its tangent vector and $\| \frac{d}{dt}
(f_t, \zeta_t) \|$ the length of its tangent vector.

\begin{theorem}\label{MainTheorem}(Main Theorem)
There is a constant $C_1 > 0$, such that:  if $(f_t, \zeta_t)$ $t_0
\leq t \leq t_1$ is a $C^1$ path in $V$, then

$$C_1 D^{3/2} \int_{t_o}^{t_1} \mu_{norm} (f_t, \zeta_t) \| \frac{d}{dt} (f_t, \zeta_t) \| dt$$

steps of projective Newton method are sufficient to continue an
approximate zero $x_0$ of $f_{t_0}$ with associated zero $\zeta_0$
to an approximate zero $x_1$ of $f_{t_1}$ with associated zero
$\zeta_1$.
\end{theorem}

 \medskip
\begin{proof}
Choose $C < \mu_0$ and $\epsilon$ small enough such that Theorem
\ref{Theorem 3} and Theorem \ref{Theorem 5} apply, $u=2C(1+\epsilon)
< \frac{3-\sqrt{7}}{4}$ and $\frac{4u}{\psi(2u)} <
\frac{1}{2(1+\epsilon)}.$ Hence, if $d (f, g) < \frac{C}{D^{1/2}
\mu_{norm}(f, \zeta)}$ ,

$d(\zeta, \eta) \le \frac{C}{D^{3/2} \mu_{norm} (f, \zeta)}$, $f
(\zeta) = 0$ ,$g(\eta)=0$

${\rm and} \ d (x, \zeta) \le \frac{C}{D^{3/2} \mu_{norm} (f,
\zeta)},$

then $d (N_g (x), \eta) \le \frac{C}{D^{3/2} \mu_{norm} (g, \eta)}$
So $N_g (x)$ is an approximate zero of $g$ with associated zero
$\eta.$

Now apply proposition 1 to produce $S_0, \dots, S_k$ and $x_0$ such
that

$$d(x_0, \zeta_{t_0})
< \frac{C}{D^{3/2} \mu_{norm} (f_{t_o}, \zeta_{t_0})}.$$

Then $x_i = N_{f_{S_i}} (x_i - 1)$ is approximate zero of $f_{S_i}$
with associated zero $\zeta_{S_i}$ and

$$
d(x_i, \zeta_{S_i}) < \frac{C}{D^{3/2} \mu_{norm} (f_{S_i},
\zeta_{S_i})}.$$

\end{proof}

\begin{corollary}
There is a constant $C_2 > 0$, such that:  if $(f_t, \zeta_t)$ $t_0
\leq t \leq t_1$ is a $C^1$ path in $V$, then

$$C_2 D^{3/2} \int_{t_o}^{t_1} \mu_{norm}^2(f_t, \zeta_t) \|\dot{f_t}\| dt$$

steps of projective Newton method are sufficient to continue an
approximate zero $x_0$ of $f_{t_0}$ with associated zero $\zeta_0$
to an approximate zero $x_1$ of $f_{t_1}$ with associated zero
$\zeta_1$.
\end{corollary}

\begin{proof}
$\|\dot{\zeta_t}\|\leq \mu_{norm}\|\dot{f_t}\|$
\end{proof}

\medskip
Theorem 4 suggests that if we wish to continue a solution $\zeta
\epsilon P (\IC^{n+1})$ of $f \epsilon P (\cH_{(d)})$ to a solution
$\eta \epsilon P (\IC^{n+1})$ of $g \epsilon P (\cH_{(d)})$ an
efficient way might be to follow a geodesic joining $(f, \zeta)$ to
$(g, \eta)$ in the metric $\|(\dot f, \dot \zeta) \|^2_k =
\mu_{norm} (f, \zeta)^2 (\| \dot f \|^2 + \| \dot \zeta \|^2)$. We
call this Riemannian metric the $\emph{condition (number)metric}$
and quickly drop the "number" from the name.

Let $\sum' \subset V = \{ (f, \zeta) \epsilon V \mid \mu_{norm} (f,
\zeta) = \infty \}$ and $W = V - \sum'$.  Note that $\mu_{norm}^2
(f, \zeta)$ is not differentiable everywhere on $W$.

\medskip

\begin{theorem}\label{Theorem 6}
$W$ is complete in the metric $\| \ \|_k$.
\end{theorem}

\medskip

\begin{lemma}
There is a constant $C > 0$ such that $(f_t, \zeta_t)$ in a $C^1$
path in $W$ , $t_0 \le t \le t_1$ of length $L$ in the $\| \ \|_k$
metric, then $\mu_{norm} (f_{t_1}, \zeta_{t_1}) \le C^{D^{3/2} L}
\mu_{norm} (f_{t_0}, \zeta_{t_0})$.
\end{lemma}

\bigskip

{\bf Proof of lemma}

From proposition 1 it follows that
$$\mu (f_{t_1}, \zeta_{t_1}) \le
(1 + \epsilon)^{(\frac{1+\epsilon}{C}) D^{3/2} \sqrt{2} L}
\mu (f_{t_0}, \zeta_{t_0})$$

For an appropriate $\epsilon, C$.  Let $C = (1 + \epsilon)^{\sqrt{2} (\frac{1+\epsilon}{C})}$

\bigskip
{\bf Proof of Theorem}

Fix $(f_0, \zeta_0)$ for example $f_{o_i} =
\frac{d_i^{1/2}}{n^{1/2}} X_i X_0^{d_i - 1} i = 1, \dots, n$

and $\zeta_0 = (1, 0, \dots, 0).$

Then $\mu_{norm} (f_0, \zeta_0) = n^{1/2}.$

Hence,
 $\mu_{norm} (f, \zeta) \le C^{D^{3/2} d \bigr( (f, \zeta_),
(f_0, \zeta_0) \bigr)}n^{1/2}.$ So

any Cauchy Sequence in $W$ stays a bounded distance away from
$\sum'$, hence in a compact region of $W$ where it converges in the
usual metric but also in the metric induced by $ \| \ \|_k$.

\end{document}